\documentclass[14pt]{extarticle}
\usepackage[cp1251]{inputenc}
\usepackage{amsmath}
\usepackage{amssymb}
\usepackage[russian,english]{babel} 

    \usepackage[pdftex]{graphicx}

\newcommand{\Z}{\mathbb{Z}}
\renewcommand{\leq}{\leqslant}
\renewcommand{\geq}{\geqslant}
\newcommand{\summ}{\sum\limits}
\newcommand{\prodd}{\prod\limits}
\newcommand{\ve}{\varepsilon}

\begin{document}
\begin{center}
\textbf{SETS AVOIDING SQUARES IN $\Z_m$ } 
\end{center}

\begin{center}
	Mikhail Gabdullin \footnote{The work is supported by the grant from Russian Science Foundation (Project 14-11-00702).}
\end{center}

\selectlanguage{english}

\begin{abstract}
We prove that for all squarefree $m$ and any set $A\subset\Z_m$ such that $A-A$ does not contain non-zero squares the bound $|A|\leq m^{1/2}(3n)^{1.5n}$ holds, where $n$ denotes the number of odd prime divisors of $m$.	 
\end{abstract}


\section{Introduction}

It was a conjecture of L.\,Lov\'asz that if $S$ is any sequence of positive integers of positive asymptotic density, then $S-S$ necessarily contains a square. 
A.~S\'ark\"ozy \cite{Sarkozy}  proved it, showing that for any $B\subset[N]$ such that $B-B$ avoids squares we have
$$ |B|\ll N(\log N)^{-1/3+\ve}.
$$
Currently the best upper bound is
$$|B|\ll \frac{N}{(\log N)^{\log\log\log\log N/12}},
$$
which was obtained by J.\,Pintz, W.~L.\,Steiger, and E.\,Szemer\'edi \cite{Pintz-Steiger-Szemeredi}. The method of that work also gives the similar upper bound to the case of $k$th powers; see \cite{Balog-Pelikan-Pintz-Szemeredi}. On the other hand, I.\,Ruzsa \cite{Ruzsa} constructed an example of a set $B\subset[N]$  which possess the mentioned property and has size $|B|\gg N^{\gamma}$, where $\gamma=\frac12(1+\frac{\log 7}{\log 65})=0,733077\ldots$.

With this connection it is natural to consider the correspondence problem in cyclic group $\Z_m$. This question is also explored by I.\,Ruzsa and M.\,Matolcsi in \cite{Ruzsa-Matolcsi}. For sets $A\subset\Z_m$ with the property that $A-A$ avoids cubic residues they showed that
$$|A|=O_{\ve}(m^{1/2+\ve})
$$
for all squarefree $m$, and
$$|A|\leq m^{1-\delta},
$$
where $\delta=0.119\dots$, for all $m$. If $A-A$ avoids squares, they proved the bounds  
\begin{equation}|A| < m^{1/2}  \label{4k+1}
\end{equation} 
for all squarefree $m$ which have prime divisors $1\pmod{4}$ only, and
$$|A|\leq me^{-c\sqrt{{\log m}}}
$$
for all squarefree $m$. 
\bigskip

In this paper we investigate the squarefree modular case for sets avoiding squares. Firstly, we would like to discuss briefly known lower bounds. It was shown by S.\,Cohen \cite{Cohen} that there exists such a set of size at least $\frac12(\log_2 m+o(1))$ for all $m$ which have prime divisors $1\pmod{4}$ only, while S.~Graham and C.\,Ringrose \cite{Graham} proved the lower bound $\log p\log\log\log p$ for infinitely many primes $m=p$.  

We present a short proof of the bound obtained in \cite{Cohen}. Let us begin with the case $m=p\equiv1 \pmod{4}$. Consider the complete graph $G=(V,E)$ with $V=\Z_p$ and the partition $E=E_1\bigsqcup E_2$, where $E_1=\{(x,y) : x-y \,\, \mbox{is a square}\}$ and $E_2=E\setminus E_1$. Then, by Ramsey's theorem for two colours (see, for instance, \cite{Tao-Vu}, Theorem 6.9), one can find a complete monochromatic subgraph $G'=(V',E')$ of our graph $G$ with $|V'|=n$ whenever $|V|=p \geq {{2n-2}\choose{n-1}}$. We thus see that there exists such a subgraph with $n\geq \frac12\log_2 p$. If $E\subset E_2$, then, obviously, the set $V'$ of all its vertices gives an example we need; if $E\subset E_1$, then for any non-residue $\xi\in\Z_p$ we get such an example in the form $\xi V'$. To get the bound for the mentioned more general case, observe that if $m=\prod_{i=1}^kp_i $ and $A_i \subset \Z_{p_i}$ possess the property that $A_i-A_i$ avoids squares, then, obviously, the set $A_1\times\ldots\times A_k$ possess it too. The claim follows.

\bigskip

Our main result is the following.

\textbf{Theorem.} 
\textit{For all squarefree $m$ and $A\subset\Z_m$ such that $A-A$ does not contain non-zero squares we have
$$ |A|\leq m^{1/2}(3n)^{1.5n},	
$$	
where $n$ denotes the number of odd prime divisors of $m$.}
\bigskip


\textbf{Corollary 1.} \textit{Let $m$ and $A$ obey the conditions of the Theorem. If $n=o(\frac{\log m}{\log\log m})$, then $$|A|\leq m^{1/2+o(1)};$$ 
if $n\leq (\frac13-\ve)\frac{\log m}{\log\log m}$, then}
$$|A|\leq m^{1-1.5\ve+o(1)}.$$ 

\textbf{Corollary 2.} \textit{We have 
$$|A|\leq m^{-c\log m/\log\log m}
$$
for all $m$ and $A$ obeying the conditions of the Theorem.}
\bigskip

The Theorem will be proven in Section 2. Corollary 1 follows immediately from the Theorem; Corollary 2 will be proven in Section 3.


\section{Proof of the Theorem}

Without loss of generality we may assume that $m$ is odd. We induct on $n$. For the case $n=1$, i.e., $m=p$ is prime, we have the bound $|A|\leq m^{1/2}$. If $p\equiv3\pmod{4}$, then $|A|\leq1$; suppose $p\equiv1\pmod{4}$. We give an elegant and folklore proof: let us assume that $|A|>m^{1/2}$ and fix a non-residue $\xi\in\Z_m$. Consider the map $\varphi\colon A^2\to\Z_p$, $\varphi(a,b)=a+\xi b$. By the pigeonhole principle, there are two distinct pairs $(a_1,b_1)$ and $(a_2,b_2)$ such that $\varphi(a_1,b_1)=\varphi(a_2,b_2)$, i.e., $\xi=(a_1-a_2)(b_2-b_1)^{-1}$, which means that at least one of the differences $a_1-a_2$ and $b_1-b_2$ is non-residue modulo $m$, and the claim follows.

 Now assume that $n\geq2$ and the claim is true for all $l<n$. 
 Let $p_1<p_2<\ldots<p_n$ be all prime divisors of $m$. Denote by $\chi_j$ quadratic character of $\mathbb{Z}_{p_j}$. Since each difference $a_1-a_2$ of distinct elements of $A$ is non-residue by at least one modulo $p_i$, we have
\begin{equation*} |A|=\summ_{a_1,a_2\in A}\prodd_{j=1}^n (1+\chi_j(a_1-a_2))
= |A|^2 + \summ_{D}\summ_{a_1,a_2\in A}\chi_D(a_1-a_2),
\end{equation*}
where $D$ runs over all non-empty subsets of $[n]=\{1,\ldots,n\}$ and $\chi_D(x)=\prodd_{j\in D}\chi_j(x)$. Denote $\sigma=1-|A|^{-1}$. Then we may rewrite the last equality as follows:
\begin{equation*}
|A|^2\sigma=-\summ_{D}\summ_{a_1,a_2\in A} \chi_D(a_1-a_2).
\end{equation*}
Using Cauchi-Schwarz, we see that
\begin{equation*} |A|^2\sigma\leq \summ_{D}|A|^{1/2}S_{D}^{1/2}, 
\end{equation*}
where
$$ S_{D}=\summ_{a\in A} \left|\summ_{b \in A} \chi_D(a-b)\right|^{2}.
$$
Thus
\begin{equation}
|A|^{3/2}\sigma\leq \summ_{D}S_{D}^{1/2}.
\label{A}
\end{equation}
Now we have to estimate the sums $S_D$. Fix a set $D$ of size $d$. Denote for the brevity
$p_D=\prodd_{j\in D}p_j$ and  
\begin{equation}
G_d=(3n)^{1.5(n-d)}. \label{G_d}
\end{equation}
For all residues $x$ modulo $p_D$ we set 
$$A_x=\{a\in A : a\equiv x \pmod{p_D}\}.
$$
One can think of elements of $A_x$ as residues modulo $mp_D^{-1}$, and the difference of distinct elements of $A_x$ is non-residue modulo $mp_D^{-1}$. Then by the induction hypothesis we have
\begin{equation*}|A_x|\leq m^{1/2}p_D^{-1/2}G_d \,. 
\end{equation*}
Obviously $A=\bigsqcup\limits_{x\in\mathbb{Z}_{p_D}} A_x$ and all elements of $A_x$ give the same contribution to $S_D$. We thus see that
\begin{multline*}
S_D = \summ_{x\in \Z_{p_D}} \summ_{a\in A_x}\left|\summ_{b \in A} \chi_D(x-b)\right|^{2} =\summ_{x\in \Z_{p_D}} |A_x|\left|\summ_{b \in A} \chi_D(x-b)\right|^{2} \leq \\
m^{1/2}p_D^{-1/2}G_d\summ_{b_1,b_2\in A}\summ_{a\in\Z_{p_D}}\prodd_{j\in D} \chi_j(a-b_1)\chi_j(a-b_2) = \\
m^{1/2}p_D^{-1/2}G_d\summ_{b_1,b_2\in A}\prodd_{j\in D}\summ_{a_j\in\Z_{p_j}} \chi_j(a_j-b_1)\chi_j(a_j-b_2). 
\end{multline*}
Let us compute the inner sum. 
For the sake of brevity we introduce the following definition: a pair $(b_1,b_2)$ is said to be \textit{special} modulo $p$ if $b_1\equiv b_2 \pmod{p}$. We have $\summ_{a\in \Z_{p_j}}\chi_j(a-b_1)\chi_j(a-b_2)=p_j-1$ if $(b_1,b_2)$ is a special pair modulo $p_j$ and 
$$\summ_{a\in \Z_{p_j}}\chi_j(a-b_1)\chi_j(a-b_2)=\summ_{a\neq b_2} \chi_j\left(1+\frac{b_2-b_1}{a-b_2}\right)=\summ_{a\neq1} \chi_j(a)=-1
$$ 
otherwise.

Denote by $B_r$ the contribution of pairs which are special exactly for $r$ modulos, $0\leq r \leq d$, to the outer sum of the bound for $S_D$ . We thus have
\begin{equation}
S_D\leq m^{1/2}p_D^{-1/2}G_d \summ_{r=0}^d B_r . \label{S_D1}
\end{equation}
Obviously,
\begin{equation}B_0\leq|A|^2. \label{B_0}
\end{equation}
To obtain an estimate for the sum $S_D$ it remains to handle with $B_r$ for $r\geq1$. Fix a set $D'\subset D$, $D'=\{i_1,\ldots,i_r\}$, of numbers of special modulus.  The contribution of pairs which are special exactly these modulus to $B_r$ is at most $p_{D'}=\prodd_{j\in D'}p_j$. The amount of such pairs does not exceed the number of solution of the congruence $x\equiv y \pmod{p_{D'}},\, x,y\in A$, which is at most $|A|m^{1/2}p_{D'}^{-1/2}G_r$ by the induction hypothesis. Thus, the contribution of pairs which are special modulus $p_j,\,j\in D'$, to $B_r$ is at most $|A|m^{1/2}p_{D'}^{1/2}G_r$. Therefore for all $r\geq1$ we have 
\begin{equation} B_r\leq |A|m^{1/2}G_r\summ_{D'\subset D, \,|D'|=r}p_{D'}^{1/2}.
\label{B_r}
\end{equation}
Substituting (\ref{B_0}) and (\ref{B_r}) into (\ref{S_D1}), we see that for all $|D|=d$ 
 $$S_D \leq m^{1/2}p_D^{-1/2}G_d|A|^2+m|A|G_d\summ_{r=1}^{d} G_r \summ_{D'\subset D,\,|D'|=r} (p_{D'}/p_D)^{1/2}, 
$$
or
\begin{equation*} S_D \leq m^{1/2}p_D^{-1/2}G_d|A|^2+m|A|G_d\summ_{r=1}^{d} G_r \summ_{D'\subset D,\,|D'|=d-r} p_{D'}^{-1/2}.
\label{S}
\end{equation*}
This implies
$$S_D^{1/2} \leq m^{1/4}p_D^{-1/4}G_d^{1/2}|A| + |A|^{1/2}m^{1/2}G_d^{1/2}\summ_{r=1}^{d}G_{r}^{1/2}\summ_{D'\subset D,\,|D'|=d-r}p_{D'}^{-1/4} . 
$$
Substituting this estimate into (\ref{A}), we obtain
\begin{equation}
|A|\sigma \leq |A|^{1/2}m^{1/4}T_1 + m^{1/2}T_2, \label{6} 
\end{equation}
where
\begin{equation*} T_1=\summ_{d=1}^{n}G_{d}^{1/2}\summ_{|D|=d}p_D^{-1/4}, \label{T_1}
\end{equation*} 
\begin{equation*}T_2=\summ_{D\subseteq[n]} G_{|D|}^{1/2}\summ_{D'\subset D}G^{1/2}_{|D|-|D'|}p_{D'}^{-1/4}. \label{T_2}
\end{equation*} 

It remains to estimate the sums $T_1$ and $T_2$. We firstly handle with $T_1$. Since $p_1\geq3$ and the function $u^{-1/4}$ is concave, we have
\begin{multline}
\summ_{j=1}^n p_j^{-1/4}\leq\summ_{j=1}^n(2j+1)^{-1/4}\leq 0.5\summ_{j=1}^n\int_{2j}^{2j+2}u^{-1/4}du=\\
\frac23((2n+2)^{3/4}-2^{3/4}) < \frac23(2n)^{3/4} < 1.13n^{3/4}.  \label{1.13}
\end{multline}
Hence, recalling the definition (\ref{G_d}) of $G_d$, 
\begin{multline}
T_1\leq \summ_{d=1}^nG_d^{1/2}\frac{1}{d!}\left(\summ_{j=1}^np_j^{-1/4}\right)^d\leq 
\summ_{d=1}^n \frac{1.13^d}{d!}(3n)^{0.75(n-d)}n^{0.75d} \\ = (3n)^{0.75n}\summ_{d=1}^{n} 3^{-0.75d}\frac{1.13^d}{d!}\leq 0.65(3n)^{0.75n}. \label{T1} 
\end{multline}
Now we are going to estimate $T_2$. We may rewrite
\begin{equation*}
T_2 = \summ_{D'\subset[n]}p_{D'}^{-1/4}\summ_{D\supset D'}G_{|D|}^{1/2}G_{|D|-|D'|}^{1/2}.  
\end{equation*}
We begin with an estimate for the inner sum. By (\ref{G_d}), we see that
\begin{multline*}
\summ_{D\supset D'}G_{|D|}^{1/2}G_{|D|-|D'|}^{1/2} =
(3n)^{1.5n}\summ_{D\supset D'} (3n)^{-1.5(|D|-|D'|/2)}\leq \\ 
(3n)^{1.5n}\summ_{r=|D'|+1}^n n^{r-|D'|}(3n)^{-1.5(r-|D'|/2)} = \\
(3n)^{1.5n+0.75|D'|}n^{-|D'|}\summ_{r=|D'|+1}^n 3^{-1.5r}n^{-r/2} \leq \\
(3n)^{1.5n+0.75|D'|}n^{-|D'|} 3^{-1.5(|D'|+1)}n^{-(|D'|+1)/2} (1-3^{-1.5}n^{-1/2})^{-1} \leq \\
0.16(3n)^{1.5n-0.75|D'|} .
 \end{multline*}
Then, thanks to (\ref{1.13}), we obtain
\begin{multline*}
T_2\leq 0.16(3n)^{1.5n}\summ_{l=0}^{n-1}(3n)^{-0.75l}\summ_{|D'|=l}p_{D'}^{-1/4} \leq \\ 0.16(3n)^{1.5n}\summ_{l=0}^{n-1} 3^{-0.75l}\frac{1.13^l}{l!} \leq 
0.27(3n)^{1.5n}.  
\end{multline*}
In light of this and (\ref{T1}), we see from (\ref{6}) that
\begin{equation*}
L:=|A|^{1/2}\left(|A|^{1/2}\sigma - 0.65m^{1/4}(3n)^{0.75n}\right)\leq 0.27m^{1/2}(3n)^{1.5n}=:R. \label{<}
\end{equation*} 
Assume that 
$$|A|> m^{1/2}(3n)^{1.5n}.
$$
But $n\geq2$; hence, $m\geq15$, $|A|\geq 6^3\sqrt{15}>100$ and $\sigma=1-|A|^{-1}\geq 0.99$. Therefore
\begin{equation*}L> (0.99-0.65)m^{1/2}(3n)^{1.5n}>R, 
\end{equation*} 
a contradiction. This completes the proof.

\section{Proof of Corollary 2}

The idea of the proof is to combine the Theorem with another upper bound on $|A|$ which is decreasing on $n$.

Denote $m'=\prod_{p|m, \newline p=3 \pmod{4}} p$. We may assume that $m' \geq m^{1/2} $ (say), since otherwise we have $|A|\leq m'(m/m')^{1/2}\leq m^{3/4}$ by (\ref{4k+1}). For similar reasons we see that it suffices to prove the claim for the case $m'=m$. 

We will use the graph theoretic approach suggested by M.\,Matolcsi and I.\,Ruzsa \cite{Ruzsa-Matolcsi}. Recall that product $(V,E)$ of directed graphs $D_i=(V_i,E_i)$, $1\leq i\leq k$, is defined as follows: we set $V=V_1\times\ldots\times V_k$ and say that an ordered pair of distinct vertices $((x_1,\ldots,x_k),(y_1,\ldots,y_k))\in V^2$ belongs to $E$ if and only if we have either $x_i=y_i$ or $(x_i,y_i)\in E_i$ for all $i$. A directed graph is called a tournament if exactly one of $(x,y)\in E$ and $(y,x)\in E$ is true for all $x\neq y$. 

We need the following result of N.\,Alon.

\textbf{Lemma} (~\cite{Alon}, Theorem 1.2) 
\textit{Let $(V_1,E_1),\ldots,(V_k,E_k)$ be directed graphs with maximum outdegrees $d_1,\ldots,d_k$ respectively and $(V,E)$ be its product. Suppose that $S$ is a subset of $V$ with the property that for every ordered pair $(u_1,\ldots,u_k)$ and $(v_1,\ldots,v_k)$ of members of $S$ we have $(u_i,v_i)\in E_i$ for some $i$. Then}
$$ |S|\leq \prod_{i=1}^k (d_i+1).	
$$	

Note that in [8] only the case $(V_1,E_1)=\ldots=(V_k,E_k)$ is considered but the proof immediately extends to different directed graphs.	For completeness, we reproduce the proof given there.

\smallskip

\textit{Proof of the lemma.} We may think of each set $V_i$ as a set of integers. Associate each member $v=(v_1,\ldots,v_k)$ of $S$ with a polynomial $P_v\in \mathbb{Q}[x_1,\ldots,x_k]$ defined by
$$P_v(x_1,\ldots,x_k)=\prod\limits_{i=1}^k\prod\limits_{j\in N(v_i)}(x_i-j),
$$
where $N(v_i)=\{u\in V_i : (v_i,u)\in E_i\}$ is the set of all out-neighbors of $v_i$.

Since $v_i\notin N(v_i)$, we see that $P_v(v_1,\ldots,v_k)\neq0$ for all $v=(v_1,\ldots,v_k)\in S$. On the other hand, by the definition of $S$, we have $P_v(u)=0$ whenever $u\in S$ and $u\neq v$. It follows that the set of polynomial $\{P_v : v\in S\}$ is linearly independent (since if $\sum_{v\in S} c_vP_v(x_1,\ldots,x_k)=0$ then, by substituting $(x_1\,\ldots,x_k)=(v_1,\ldots,v_k)$ we conclude that $c_v=0$). But each $P_v$ is a polynomial of degree at most $d_i$ in variable $x_i$; hence, the number of these polynomials does not exceed the dimension of the space of polynomials in $k$ variables with this property, which is $\prod_{i=1}^k(d_i+1)$. This concludes the proof.

\bigskip 

Now assume that $A\subset\Z_m$ is such that $A-A$ does not contain non-zero squares. We consider the product $(\Z_m,E)$ of the tournaments $(\Z_p, E_p)$, $p|m$, where $(x,y)\in E_p$ iff $x-y$ is a square in $\Z_p$ (recall that we assume all $p$ to be $3\pmod{4}$).  Then for any $a,b$ we can find $p|m$ with $(a-b)\pmod{p}\in E_p $ (since $(b,a)\notin E$). We thus see from the lemma that $|A|\leq \prod_{i=1}^n (p_i+1)/2=m2^{-n}\prod_{i=1}^n(1+1/p_i)\leq m2^{-cn}$ for some $c>0$. Combining this with the Theorem, we get $|A|\leq m\cdot\min(2^{-cn},m^{-1/2}(3n)^{1.5n})$, and the claim follows.

\end{document}